\documentclass[english, aps, showpacs, preprint, nobibnotes]{revtex4-1}

\usepackage{amsmath,amssymb}

\usepackage[utf8x]{inputenc}

\usepackage[T1]{fontenc}
\usepackage{float}
\usepackage{mathtools}
\usepackage{amsmath}
\usepackage{amssymb}
\usepackage{graphicx}
\usepackage{nicefrac}
\usepackage{resizegather}
\usepackage{breakcites}

\begin{document}
	
	\title{Emergence of Subcritical Bifurcations in a System of Randomly Coupled Supercritical Andronov-Hopf Oscillators: A Potential Mechanism for Neural Network Type Switching}
	\date{\today}
	
	\author{Keith Hayton}
	\email{khayt86@gmail.com }
	\thanks{author contributed equally}
	\author{Dimitrios Moirogiannis}
	\email{dmoirogi@gmail.com}
	\thanks{author contributed equally}
	\affiliation{Center for Studies in Physics and Biology, The Rockefeller University}

	\begin{abstract}
		
	Experimental evidence suggests that the computational state of cortical systems change according to behavioral and stimulus context. However, it is still unknown what mechanisms underlie this adaptive processing in cortical circuitry. In this paper, we present a model of randomly coupled supercritical Andronov-Hopf oscillators which can act as an adaptive processor by exhibiting drastically different dynamics depending on the value of a single network parameter. Despite being only composed of supercritical subunits, the full system can exhibit either supercritical or subcritical Andronov-Hopf bifurcations. This model might provide a novel mechanism for switching between globally asymptotically stable and nonhyperbolic neural network types in pattern recognition theory.
		
	\end{abstract}

	\maketitle
	
	\section*{Introduction}
	
	A growing body of evidence suggests that cortical areas are not fixed in function but instead act as adaptive processors \cite{kirst2017} with intrinsic cortical circuits that are strongly modulated by feedback from higher cortical areas \cite{gilbert2007, gilbert2013} and incoming sensory information \cite{kapadia1999dynamics, sceniak1999contrast, nauhaus2009stimulus, yan2012input, Hayton2018a}. This is in contrast to the classical feedforward view of the cortex as a hierarchical series of cortical areas which encode an increasingly more complex description of sensory input \cite{gilbert2007, gilbert2013}. How the brain is able to dynamically reconfigure it's computational state is largely unknown and a matter of active research \cite{kirst2017}.
	
	In this paper, we present a toy model of cortical circuitry consisting of a network of damped, identical oscillators coupled through a \textit {Wigner real symmetric matrix} \cite{Tao2008circular, Anderson2010, Tao2012Random, O'Rourke2016Eigenvectors}. Each oscillator in the network is poised at a supercritical \textit{Andronov-Hopf bifurcation} \cite{poincare1893, Hopf1942, Andronov2013}. To study the dynamics \cite{Ioos1999, wiggins2003introduction, Kuznetsov2013} of the system, we follow the statistical \textit{center manifold reduction} approach introduced in \cite{Moirogiannis2019}. We show that by changing a single network parameter the full system can undergo either a supercritical or subcritical Andronov-Hopf bifurcation. The transformation between supercritical and subcritical occurs despite the system being exclusively composed of supercritical oscillators. This model is an example of how a network of fixed computational subunits can act as an adaptive processor by drastically switching it's dynamic regime without affecting the subunit dynamics. Additionally, we also explore how this network model might also act as a neuroscience inspired adaptive pattern recognition processor by switching between neural network types, specifically \textit{globally asymptotically stable} and \textit{nonhyperbolic}. 
	
	\section*{Analysis}
	
	We consider the following system of $N$ coupled, damped, identical oscillators:
	
	\begin{eqnarray}
	\dot{x}_{i}&=&y_{i}+ax_{i}^2+bx_{i}^3+\underset{j}{\sum}M_{ij}x_{j}\label{eq:system}\\
	\dot{y}_{i}&=&-x_{i}\nonumber 
	\end{eqnarray}
	
	\noindent , where $i\in\left\{ 1,\ldots, N\right\} $ and $a,b \in \mathbb{R}$ are system parameters. The individual subunits
	
	\begin{eqnarray}
	\dot{x}_{i}&=&{y}_{i}+a{x}_{i}^2+b{x}_{i}^3\label{eq:units}\\
	\dot{y}_{i}&=&-x_{i}\nonumber 
	\end{eqnarray}	
	
	\noindent are linearly coupled through the variables $x_{i}$, with the connectivity matrix $M \in \mathbb{R}^{N \times N}$ given by a random matrix with all eigenvalues negative and contained within a compact set in the complex plane except for a single leading eigenvalue $\lambda$, which will act as a bifurcation parameter. The construction of such a matrix is outlined in \cite{Moirogiannis2019}. 
    
    In general, a system 
    
    \begin{eqnarray}
    	\dot{x}&=& y+f\left(x,y\right)\label{eq:Hopf}\\
    	\dot{y}&=&-x+g\left(x,y\right)\nonumber
    \end{eqnarray}
    	
    \noindent is poised at an Andronov-Hopf bifurcation if $f$ and $g$, which capture all nonlinearities of order $\geq2$, ensure that the first Lyaponuv coefficient $l_1(0)=\dfrac{1}{16}\left(f_{xxx}+f_{xyy}+g_{xxy}+g_{yyy}\right)+\dfrac{1}{16}\left(f_{xy}\left(f_{xx}+f_{yy}\right)-g_{xy}\left(g_{xx}+g_{yy}\right)-f_{xx}g_{xx}+f_{yy}g_{yy}\right) \neq 0$ \cite{{Kuznetsov2013}}. The bifurcation is subcritical (supercritical) iff $l_1(0) > 0$ ($l_1(0) < 0$).
    
    For the individual subunits described in (\ref{eq:units}), the first Lyapunov coefficient is given by $l_1(0)=\dfrac{1}{16}\left(f_{xxx}+f_{xyy}\right)+\dfrac{1}{16}\left(f_{xy}\left(f_{xx}+f_{yy}\right)\right)=\dfrac{1}{2}b$, where $b$ is the coefficient of the cubic term in (\ref{eq:units}). Therefore, the type of bifurcation that a single subunit is poised at is determined exclusively by $b$, with the bifurcation being subcritical (supercritical) iff $b>0$ ($b<0$); the coefficient $a$ of $x^2$ in (\ref{eq:units}) does not play a role.
    
	We will now proceed to calculate the Lyapunov coefficient for the full system. Using the statistical center manifold reduction technique outlined in \cite{Moirogiannis2019}, it has been shown that at the bifurcation point $\lambda=0$, the system of equations in (\ref{eq:system}) can be reduced to a set of restricted equations on the center manifold. Up to third order, these reduced equations are given by:
	
	\begin{eqnarray}
		\nonumber\dot{x}&=&y+a\sqrt{N}\Gamma_{1}^{2\delta_{1,\bullet}}x^{2}+\left(b{\sqrt{N}}^{2}\Gamma_{1}^{3\delta_{1,\bullet}}+2a{\sqrt{N}}^{2}\underset{k=2}{\overset{N}{\sum}}\Gamma_{1}^{\delta_{1,\bullet}+\delta_{k,\bullet}}\alpha_{0}^{k}\right)x^{3}+\\ 
		 +&&\left(2a{\sqrt{N}}^{2}\underset{k=2}{\overset{N}{\sum}}\Gamma_{1}^{\delta_{1,\bullet}+\delta_{k,\bullet}}\alpha_{1}^{k}\right)x^{2}y+\left(2a{\sqrt{N}}^{2}\underset{k=2}{\overset{N}{\sum}}\Gamma_{1}^{\delta_{1,\bullet}+\delta_{k,\bullet}}\alpha_{2}^{k}\right)xy^{2}+O\left(4\right)\\ 
		\nonumber \dot{y}&=&-x
	\end{eqnarray}
	
	\noindent , where the $\Gamma$'s are random variables which depend on the structure of the eigenvectors and are defined in \cite{Moirogiannis2019}. In the limit $N \rightarrow \infty$, the $\Gamma$'s can be evaluated using moments of the statistical distribution of the eigenvectors \cite{O'Rourke2016Eigenvectors, Moirogiannis2019}.

	By analyzing the reduced equations, it can be proved that the system undergoes an Andronov-Hopf bifurcation when $\lambda=0$. Observe that the reduced equations are of the form (\ref{eq:Hopf}), and thus, the Lyapunov coefficient can be easily calculated: $l_1(0)=b\dfrac{3}{8}{\sqrt{N}}^{2}\Gamma_{1}^{3\delta_{1,\bullet}}+a^{2}{\sqrt{N}}^{2}\underset{k=2}{\overset{N}{\sum}}\Gamma_{1}^{\delta_{1,\bullet}+\delta_{k,\bullet}}\Gamma_{k}^{2\delta_{1,\bullet}}\dfrac{-\lambda_{k}}{4\lambda_{k}^{2}+9}\neq0$. As $l_1(0)$ only vanishes along a single parabolic curve in the $a,b$ parameter plane, this implies that the system typically undergoes a Andronov-Hopf bifurcation as $\lambda$ crosses the imaginary axis. The bifurcation is subcritical (supercritical) iff $l_1(0)>0$ ($l_1(0)<0$). 
	
	Moreover, assuming the connection matrix $M$ is normal, it can be easily shown that the term ${\sqrt{N}}^{2}\underset{k=2}{\overset{N}{\sum}}\Gamma_{1}^{\delta_{1,\bullet}+\delta_{k,\bullet}}\Gamma_{k}^{2\delta_{1,\bullet}}\dfrac{-\lambda_{k}}{4\lambda_{k}^{2}+9}$ is typically positive. First, since $M$ is normal, the $\Gamma$'s can be greatly simplified: $\Gamma_{1}^{\delta_{1,\bullet}+\delta_{k,\bullet}}=\Gamma_{k}^{2\delta_{1,\bullet}}=\underset{\phi}{\sum}V_{\phi,k}V_{\phi,1}^2$ where $(V_{1,k},\dots ,V_{N,k})^t$ is the $k$th vector in the orthonormal base of eigenvectors (with eigenvalue $\lambda_{k}$) \cite{Moirogiannis2019}. It then follows that ${\sqrt{N}}^{2}\underset{k=2}{\overset{N}{\sum}}\Gamma_{1}^{\delta_{1,\bullet}+\delta_{k,\bullet}}\Gamma_{k}^{2\delta_{1,\bullet}}\dfrac{-\lambda_{k}}{4\lambda_{k}^{2}+9}={\sqrt{N}}^{2}\underset{k=2}{\overset{N}{\sum}}(\underset{\phi}{\sum}V_{\phi,k}V_{\phi,1}^2)^2\dfrac{-\lambda_{k}}{4\lambda_{k}^{2}+9}> 0$ almost surely. We also have that the other term in the Lyapunov coefficient is positive $\Gamma_{k}^{3\delta_{1,\bullet}}>0$, since $\Gamma_{k}^{3\delta_{1,\bullet}}=\underset{\phi}{\sum}V_{\phi,k}V_{\phi,1}^3=\underset{\phi}{\sum}V_{\phi,1}^4>0$.
	
	We now examine how the bifurcation of the full system and individual subunits vary with parameters $a$ and $b$. If $b>0$, the full system along with each subunit ($l_1(0)_{subunit}=\dfrac{1}{2}b$) is subcritical regardless of the value of $a$. When $b<0$, the subunits are always supercritical, but the full system can be either subcritical or supercritical depending on the value of $||a||$: the full system is subcritical (supercritical) iff $||a||>\gamma \sqrt{-b}$ ($||a||<\gamma \sqrt{-b}$), where $\gamma=\sqrt{\nicefrac{\frac{3}{8}\underset{\phi}{\sum}V_{\phi,1}^4} {\left(\underset{k=2}{\overset{N}{\sum}}(\underset{\phi}{\sum}V_{\phi,k}V_{\phi,1}^2)^2\frac{-\lambda_{k}}{4\lambda_{k}^{2}+9}\right)}}$ . Thus, it is possible for the full system to undergo a subcritical Andronov-Hopf bifurcation despite each individual subunit being supercritical.  
	
	\section*{Adaptive Pattern Recognition}
	
	We now explain how our results might be applicable to pattern recognition and switching of neural network types. 	Assuming static input patterns, neural networks can be broadly categorized into three types: \textit{multiple attractor} (MA), \textit{globally asymptotically stable} (GAS), and \textit{nonhyperbolic} (NH) neural networks \cite{hoppensteadt2012weakly}. MA-type networks, which include Hopfield and Cohen-Grossberg networks, are characterized by each input pattern converging to one of many multiple attractors; each attractor corresponds to a different stored memory pattern. GAS-type networks, on the other hand, store memories by associating a given input pattern with the globally asymptotically stable state of the system. In this case, a memory is stored as a specific location of the attractor, and an input pattern is presented to the network as a parameter which controls the attractor location \cite{hoppensteadt2012weakly}. 
	
	Finally, NH-types have nonhyperbolic equilibriums which affect the global dynamics of the system. Input in NH types is given to the system in the form of a bifurcation parameter that perturbs the nonhyperbolic equilibrium, causing it to lose stability or disappear altogether. Only a local analysis near the equilibrium is needed to understand global behavior, since the equilibrium is a point of convergence of the curves, called separatrices, separating the different attraction domains of the system \cite{hoppensteadt2012weakly}. Nonhyperbolic equilibrium dynamics have been the focus of an extensive number of studies in neuroscience. These studies include entire hemisphere ECoG recordings \cite{Solovey2012, Alonso2014, solovey2015loss}, experimental studies in premotor and motor cortex \cite{churchland2012}, theoretical \cite{seung1998continuous} and experimental studies \cite{seung2000stability} of \textit{slow manifolds} (a specific case of center manifolds) in oculomotor control, slow manifolds in decision making \cite{machens2005}, Andronov-Hopf bifurcation \cite{poincare1893, Hopf1942, Andronov2013} in the olfactory system \cite{freeman2005metastability} and cochlea \cite{choe1998model, eguiluz2000essential, camalet2000auditory, kern2003essential, Duke2003, Magnasco2003, Hayton2018b}, a nonhyperbolic model of primary visual cortex (V1) \cite{yan2012input, Hayton2018a}, and theoretical work on regulated criticality \cite{bienenstock1998regulated}.  
	
 	Following the lead of \cite{hoppensteadt2012weakly}, we consider a specific case of (\ref{eq:system}), where the $N$ damped oscillators in the system, given by (\ref{eq:units}), are neural relaxation oscillators with input matrix $R=\begin{pmatrix}r_{1} & & 0 \\ & \ddots\\ 0 & & r_{N} \end{pmatrix}$, where $r_{i}$ is the external input to the $i$-th oscillator, and $M=R+C$, where $C$ represents the linear connectivity between the $x_{i}$'s. In the case of neural relaxation oscillators, $x_{i}$ and $y_{i}$ represent rescaled populations of excitatory and inhibitory populations, respectively. As in \cite{hoppensteadt2012weakly}, we assume Dales principle holds, $C_{ij}{\geq0}$ ${\forall}i,j$, and because of Hebbian learning mechanisms, $C$ is symmetric.
	
	It has been shown, using the Perron-Frobenius theorem, that for any input matrix $R$, the leading eigenvalue $\lambda$ of $R+C$ is simple, and if $\lambda=0$, the full system undergoes an Andronov-Hopf bifurcation \cite{hoppensteadt2012weakly}. Using this result, it has also been shown that such a system has the potential to perform pattern recognition on a set of input vectors. The exact construction of the pattern recognition algorithm is still an open problem as it requires the construction of a specific matrix for all input vectors, however, once this matrix is constructed, the resulting algorithm will depend heavily on the type of bifurcation the neural network undergoes. If the bifurcation is supercritical, each input vector is mapped to a stable limit cycle with a location corresponding to one of the memorized pattern vectors. This is an example of a GAS-type neural network. On the other hand, if the bifurcation is subcritical, the trajectory corresponding to the initial input vector leaves the vicinity of the equilibrium point along a direction associated with one of the memorized patterns, which is an example of a NH-type neural network.
	
	Our earlier analysis of the coupled supercritical oscillators with normal connectivity matrix directly applies to this setting, since according to the Perron-Frobenuis theorem, a change in input corresponds to moving the leading eigenvalue as we assumed above in our calculations. Thus, the network parameter $||a||$, can switch the neural network type between GAS (supercritical) and NH (subcritical). 
	
	\section*{Conclusion}
	
	Utilizing center manifold reduction theory and Lyapunov coefficients, we have shown that a network of randomly coupled supercritical Andronov-Hopf oscillators can undergo either a supercritical or, surprisingly, a subcritical Andronov-Hopf bifurcation. We have also discussed how this switching between supercritical and subcritical regimes could potentially lead to a mechanism for neural network type switching in cortex. It is feasible that feedback connections from higher cortical areas could alter the activity of inhibitory interneurons in the network \cite{gilbert2007, gilbert2013}, leading to global changes in network parameters such as $||a||$. More experiments are needed to determine if the switching of dynamical regimes as presented in this paper might underlie any of the adaptive processing behaviors observed in intrinsic cortical circuitry.

\end{document}